\documentclass[11pt]{amsart}
\usepackage{amssymb}

\numberwithin{equation}{section}

\newtheorem{theorem}{Theorem}
\newtheorem{Theorem}[theorem]{Theorem}

\newtheorem{proposition}[theorem]{Proposition}
\newtheorem{lemma}[theorem]{Lemma}

\theoremstyle{definition}
\newtheorem{definition}[theorem]{Definition}
\theoremstyle{remark}
\newtheorem{remark}[theorem]{Remark}

\newcommand{\cL}{\mathcal{L}}

\newcommand{\calS} {\mathcal{S}}
\newcommand{\cS} {\mathcal{S}}
\newcommand{\cT}{\mathcal{T}}

\newcommand{\R}{\mathbb{R}}
\newcommand{\C}{\mathbb{C}}

\newcommand{\fg}{\mathfrak{g}}
\newcommand{\fm}{\mathfrak{m}}

\newcommand{\fk}{\mathfrak{k}}
\newcommand{\ft}{\mathfrak{t}}

\newcommand{\fu}{\mathfrak{u}}
\newcommand{\fsp}{\mathfrak{sp}}

\newcommand{\SP}{\operatorname{Sp}}

\newcommand{\la}{\langle}
\newcommand{\ra}{\rangle}
\newcommand{\inv}{^{-1}}
\renewcommand{\ker}{ \operatorname {ker}}

\begin{document}

\title{ Normal modes in symplectic stratified spaces}
\author{Eugene Lerman }

\address{Department of
Mathematics, University of Illinois, Urbana, IL 61801}
\email{lerman@math.uiuc.edu}

\thanks{Supported by the NSF grant DMS-9803051}

\date{\today}

\begin{abstract}
We generalize the Weinstein-Moser theorem on the existence of
nonlinear normal modes (i.e., periodic orbits) near an equilibrium in
a Hamiltonian system to a theorem on the existence of relative
periodic orbits near a relative equilibrium in a Hamiltonian system
with continuous symmetries.

More specifically we significantly improve a result proved earlier
jointly with Tokieda: we remove a  strong technical hypothesis
on  the symmetry group.	
\end{abstract}

\maketitle
\section{Introduction}
The goal of the paper is to generalize the Weinstein-Moser theorem
(\cite{W1, Ms, W2, MnRS, Ba} and references therein) on the nonlinear
normal modes (i.e., periodic orbits) near an equilibrium in a
Hamiltonian system to a theorem on the existence of relative periodic
orbits (the  normal modes of the title) near relative
equilibria of a symmetric Hamiltonian system.

More specifically let $(M, \omega _M)$ be a symplectic manifold with a proper
Hamiltonian action of a  Lie group $G$ and a corresponding
moment map $\Phi: M \to \fg^{\ast}$.  Assume that the moment map is
equivariant.  Let $h\in C^{\infty}(M)^G$ be a $G$-invariant
Hamiltonian.  We will refer to the tuple $(M, \omega _M , \Phi: M \to
\fg^*, h \in C^\infty (M)^G)$ as a {\bf symmetric Hamiltonian system}.

The main result of the paper is the following theorem (the terms used
in the statement are explained below):

\begin{theorem}
	\label{main_theorem} 
Let $(M, \omega _M , \Phi: M \to \fg^*, h \in C^\infty (M)^G)$ be a
symmetric Hamiltonian system.  Suppose $m\in M$ is a relative
equilibrium for the system such that the the coadjoint orbit through
$\mu = \Phi (m)$ is locally closed.  Suppose further that there exists
a symplectic slice $ \Sigma \hookrightarrow M$ through $m$ such that
the Hessian $d^2 (h|_\Sigma) (m)$ is positive definite.

Then for every sufficiently small $E>0$ the set $\{h = E + h (m) \}
\cap \Phi\inv (\mu)$ (if nonempty) contains a relative periodic
orbit of $h$.  

More precisely, let $h_\mu$ denotes the reduced Hamiltonian on the
reduced space $M_\mu = \Phi\inv (G\cdot \mu)/G$.  Then for every
sufficiently small $E>0$ the set $\{h_\mu = E + h (m)\}$ contains $N$
or more periodic orbits of the reduced Hamiltonian $h_\mu$, where $N$
is the Liusternik-Schnirelman category of the union of closed strata
in the symplectic link of the point $\bar{m} \in M_\mu$ corresponding
to $m\in M$.
\end{theorem}

Recall that given a symmetric  Hamiltonian system $(M, \omega _M ,
\Phi: M \to \fg^*, h \in C^\infty (M)^G)$, a point $m\in M$ is a
{\bf relative equilibrium} of the Hamiltonian vector field
$X_h$ of $h$ (a relative equilibrium of $h$ for short), if the
trajectory of $X_h$ through $m$ lies on the $G$ orbit through $m$.
Equivalently, since the flow of $X_h$ is $G$ equivariant, it descends to
a flow on the quotient $M/G$; $m$ is a {\bf relative equilibrium} if
the corresponding point $\bar{m} \in M/G$ is fixed by the induced
flow.  Thus if $m$ is a relative equilibrium, then the whole $G$ orbit
$G\cdot m$ consists of relative equilibria.  Similarly, we say that a
trajectory $\gamma (t)$ a $G$ invariant Hamiltonian vector field $X_h$
is a {\bf relative periodic orbit} (r.p.o.) if there is $T> 0$ and $g\in G$
such that $\gamma (T) = g\cdot \gamma (0)$. Equivalently $\gamma (t)$
is relatively periodic if the corresponding trajectory $\bar{\gamma}
(t)$ in $M/G$ is periodic.

Recall next that if a Lie group $G$ acts properly on a manifold $M$,
then at every point $x\in M$ there is a {\bf slice} for the action of
$G$, i.e., there is a submanifold $S$ passing through $x$, which is
invariant under the action of the isotropy group $G_x$ of $x$, is
transverse to the orbit $G\cdot x$ and such that the open set $G\cdot S$
is diffeomorphic to the associated bundle $G\times _{G_x} S$.

If additionally the manifold $M$ has a $G$ invariant symplectic form
$\omega _M$, then the local normal form theorem of Marle and of
Guillemin and Sternberg \cite{marle:model, GS} guaranties that we can
find a symplectic submanifold $\Sigma$ passing through $x$ which is
$G_x$ invariant with the property that the tangent space $T_x
\Sigma$ is a maximal symplectic subspace of the tangent space to the
slice $T_x S$.  Moreover, $\Sigma $ can be chosen to be $G_x$
equivariantly symplectomorphic to a ball about the origin in a linear
symplectic representation of $G_x$ on $T_x \Sigma$.  Such a
submanifold is called a {\bf symplectic slice}.  

The symplectic link of a point $\bar{m}$ in a reduced space $M_\mu$ is
a symplectic stratified space which is an invariant of a singularity
of $M_\mu$ at $\bar{m}$.  A precise definition is given later.  If the
space $M_\mu$ is smooth near $\bar{m}$ then the symplectic link of
$\bar{m}$ is smooth: it is the complex projective space $\C P^{n-1}$
where $n = \frac{1}{2} \dim M_\mu$.

If the isotropy group of the point $m$ is trivial, then the reduced
space $M_\mu$ is smooth near the corresponding point $\bar{m}$.
Moreover the symplectic slice $\Sigma$ through $m$ is symplectomorphic
to a neighborhood of $\bar{m}$ in $M_\mu$, and under the
identification of $\Sigma$ with an open subset of the reduced space
the restriction $h|_\Sigma$ is the reduced Hamiltonian $h_\mu$.  In
this case Theorem~\ref{main_theorem} reduces to a theorem of Weinstein
on the existence of nonlinear normal modes of a Hamiltonian system
near an equilibrium:

\begin{theorem}[Weinstein, \cite{W1}]\label{theorem_W}
Let $h$ be a Hamiltonian on a symplectic vector space $V$ such that
the differential of $h$ at the origin $dh(0)$ is zero and the Hessian
at the origin $d^2h(0)$ is positive definite.  Then for every small
$\varepsilon > 0$, the energy level $h\inv(h(0) + \varepsilon)$
carries at least $\frac{1}{2}\dim V$ periodic orbits of (the
Hamiltonian vector field of) $h$.
\end{theorem}

On the other hand, if $m$ is a {\it singular\/} point of the moment
map $\Phi$, then the reduced space $M_\mu$ at $\mu$ is a stratified
space, and the reduced dynamics preserves the stratification
\cite{AMM, SL, BL}.  Unless the stratum through $\bar{m}$ is an
isolated point, we have again $\frac{1}{2}\dim (\text{stratum})$
families of r.p.o.'s, provided appropriate conditions hold on the
Hessian of the restriction of the reduced Hamiltonian to the stratum.
But what if the stratum through $\bar{m}$ is a  point?

Recently Tokieda and I showed that in the singular case as in the
regular case there are relative periodic orbits near the relative
equilibrium provided a certain quadratic form is definite and the
isotropy group $G_{\mu}$ of $\mu$ is a torus \cite{LT}. Here as above
$\mu$ is the value of the moment map on the relative equilibrium. The
proof amounted to a reduction to the case where the full symmetry
group is a torus followed by computation in `good coordinates.'  The
computation allowed us to reduce our problem to Weinstein's nonlinear
normal mode theorem.\footnote{The assumption that $G_\mu$ is a torus
is not essential.  For the proof to work it is sufficient to assume
that $\mu$ is split in the sense of \cite{GLS} and that $G_\mu$ splits
up to a finite cover: $G_\mu = (G_m \times H)/\Gamma$ where $G_m$ is
the isotropy group of the relative equilibrium, $H$ a complementary
subgroup and $\Gamma$ is a finite group.} This left open a question:
\begin{itemize}
\item  Can  the assumption on the isotropy group of the value
of the moment map at the relative equilibrium be removed?
\end{itemize}
Theorem~\ref{main_theorem} answers the question affirmatively. The
guiding principle comes from \cite{SL}: a symplectic quotient is
locally modeled on a symplectic quotient of a symplectic slice
(cf. Theorem~15 in \cite{BL}).

We end the paper with Proposition~\ref{last prop} which provides
a practicable method for checking the existence of a symplectic slice
$\Sigma$ through a relative equilibrium $m$ so that the Hessian $d^2
(h|_\Sigma) (m)$ is positive definite.

\subsection*{Acknowledgments}  I thank Chris Woodward for a very
useful discussion and Yuli Rudiak for patiently answering 
my e-mail messages.

As I was writing up this paper I discovered that Ortega and Ratiu
have independently obtained a similar result \cite{OR}.

\section{ Proof of Theorem~\ref{main_theorem}}

Our first step is to reduce the proof to a special case where the
manifold $M$ is a symplectic vector space and the relative equilibrium
is an equilibrium.   Let $(M,\omega _M , \Phi: M \to \fg^*, h \in
C^\infty (M)^G)$ be a symmetric Hamiltonian system and suppose $m\in
M$ is a relative equilibrium for the system.  Then the restriction
$h_\Sigma$ of $h$ to the symplectic slice $\Sigma$ through $m$
satisfies $d h_\Sigma (m) = 0$.  This is because the Hamiltonian
vector field of $h$ points along the the group orbit $G\cdot m$ and
the tangent space to the orbit $T_m (G\cdot m)$ lies in the symplectic
perpendicular to the symplectic slice directions.  Consequently the
Hessian $d^2 h_\Sigma (m) = d^2 (h|_\Sigma)(m)$ is well-defined.

Recall next that if the coadjoint orbit through $\mu = \Phi (m)$ is
locally closed, then the reduced space $M_\mu := \Phi \inv (G\cdot
\mu) /G$ is locally isomorphic, as a symplectic stratified space, to
the reduction at zero of a symplectic slice $\Sigma$ through $m$ by
the action of the isotropy group $G_m$: see, for example
\cite[Theorem~15]{BL}. Moreover, it follows from the proof of
Theorem~15, {\it op.~cit.}, that if $h$ is any $G$ invariant
Hamiltonian on $M$ then the corresponding reduced Hamiltonian $h_\mu$
on $M_\mu$ near $\bar{m} = G\cdot m$ can also be obtained by first
restricting $h$ to $\Sigma$ and then carrying out the reduction by the
group $G_m$.

Since the symplectic slice is equivariantly symplectomorphic to a ball
in a symplectic vector space with a linear symplectic action of a
compact Lie group, it follows that in order to prove
Theorem~\ref{main_theorem}, it suffices to prove the following special
case:

\begin{Theorem}
	\label{theorem_linear} 
Let $(V, \omega _V)$ be a symplectic vector space with a linear
symplectic action of a compact Lie group $K$, and let $\Phi_V :V \to
\fk^*$ denotes the corresponding homogeneous moment map.   Assume that
$\Phi_V \inv (0) \smallsetminus \{0\}$ is nonempty, i.e., assume that
the reduced space $V_0 =\Phi_V \inv (0)/K$ is not one point.  Suppose
$h_V \in C^\infty (V)^K$ is a $K$ invariant Hamiltonian with $d h_V
(0) = 0$, and suppose the Hessian $d^2 h_V (0)$ is positive definite.
Then for every $E > h_V (0) $ sufficiently close to $h_V (0)$ the set
$\{h_0 = E\}$ in the reduced space $V_0$ contains $N$ or more periodic
orbit of the reduced Hamiltonian $h_0$. Here $N$ is the
Liusternik-Schnirelman category of the union of closed strata in the
symplectic link of the point $* \in V_0$ corresponding to $0\in V$.
\end{Theorem}

The idea of the proof of Theorem~\ref{theorem_linear} is
straightforward.  Consider the quadratic part $q(x)$ of the
Hamiltonian $h_V$ at zero, that is, let $q(x) = d^2 h(0)(x, x)$.  Let
$q_0$ denote the corresponding reduced Hamiltonian on reduced space
$V_0$.  We will see that for every sufficiently small $E>0$ and for
certain  strata $\cT$ of $V_0$ the manifolds  
$$
\{q_0 = E\} \cap \cT
$$
contain weakly nondegenerate periodic manifolds $C\subset \{q_0 =
E\} \cap \cT$ of $q_0$.  Then by a theorem of Weinstein
\cite[p. 247]{W2}, every such {\em compact}  manifold $C$ gives rise to
Cat$(C/S^1)$ periodic orbits of the reduced Hamiltonian $h_0$ in the
stratum $\cT$, where Cat denotes the Liusternik-Schnirelman category.

Recall a characterization of weakly nondegenerate periodic
manifolds given in a corollary on p.~246 of \cite{W2}, which we
take as our definition.
\begin{definition}
Let $(N, \omega_N, h)$ be a Hamiltonian system.  A submanifold $C$ of
$N$ consisting of periodic orbits of the Hamiltonian vector field
$X_h$ of $h$ is {\bf weakly nondegenerate} iff for each orbit $c$ in
$C$
\begin{itemize}
\item[(i)] $X_h (c) \not = 0$ and

\item[(ii)] the space $\{ x \in T_{c(0)} (h\inv (E)) \mid \, x - Px$ is
a multiple of $X_h (c(0))\}$ has the same dimension as $C$.  Here
$E = h(c(0))$, and $P: T_{c(0)} (h\inv (E))\to T_{c(0)} (h\inv (E)) $
denotes the linearization of the Poincar\'e map along the periodic
orbit $c$. 
\end{itemize}
\end{definition}
Thus to prove Theorem~\ref{theorem_linear} it is enough to establish
the existence of compact weakly nondegenerate periodic manifolds of
the reduced Hessian $q_0$ and to estimate the Liusternik-Schnirelman
category of the quotients of these manifolds by $S^1$.  

We now proceed with a proof of Theorem~\ref{theorem_linear}.  Since
the function $q$ is quadratic, its Hamiltonian vector field $X_q$ is
linear, hence of the form $X_q (x) = \xi x$ for some linear map $\xi
\in \fsp (V, \omega)$, the Lie algebra of the symplectic group $\SP
(V,\omega)$.   Since $q$ is definite, $\xi$ must lie in a compact Lie
subalgebra of $\fsp (V, \omega)$; in particular  the closure of $\{ \exp t\xi
\mid \, t\in \R\}$ is a torus $T \subset \SP (V, \omega)$.   Since $q$ is
$K$-invariant, the groups $K$ and $T$ commute in $\SP (V, \omega)$.
Since both groups are compact, we may assume that $V= \C ^n$ ($n=
\frac{1}{2} \dim V$), that $K$ is a subgroup of $U(n)$ and that $T$ is
contained in the standard maximal torus of $U(n)$.

\begin{lemma}\label{lemma1}
Let $(V, \omega, K, \Phi: V \to \fk^*, q(x) \in C^\infty (V)^K)$ be as
above.  Then the Hamiltonian $q$ has no relative equilibria in the set
$\Phi\inv (0)\smallsetminus \{0\}$.  Consequently the action of the
torus $T$ generated by $q$ on the reduced space $V_0$ has no fixed
points in $V_0 \smallsetminus \{*\} = (\Phi\inv (0)\smallsetminus
\{0\})/K$.
\end{lemma}
\begin{proof}
If $v\in \Phi\inv (0)\smallsetminus \{0\}$ is a relative equilibrium
then 
$$ 
	d \la \Phi, \eta \ra (v) = dq (v) 
$$ 
for some $\eta \in \fk$.  Since $\Phi$ is quadratic homogeneous, we
have $tv \in \Phi\inv (0)$ for all $t>0$.  Hence $v \in \ker d \la
\Phi, \eta \ra (v)$. On the other hand, since $q$ is definite, the ray
$\{tv \mid t> 0\}$ is transverse to the level set $\{q = q (v)\}$,
hence $v\not \in \ker dq (v)$.  Contradiction.
\end{proof}

\subsection*{The structure of the symplectic quotient $V_0$.}
Next we tersely recall a number of results explained in \cite{SL}.
Suppose as above that $T$ is a subtorus of the maximal torus of $U(n)$
and that $K$ is a closed subgroup of $U(n)$ which commutes with
$T$. Let $U$ denote the central circle subgroup of $U(n)$.  Then the
actions of $U$ and $K$ on $\C ^n$ commute, the action of $U$ is
Hamiltonian and a moment map $f$ for the action of $U$ on $\C^n$ can
be taken to be $f(z) = ||z||^2$.

Since by assumption $\Phi\inv (0)\smallsetminus \{0\} \not = \emptyset$,
the group $U$ is not contained in $K$.  In fact the Lie algebras of $U$
and $K$ intersect trivially.  Let us prove this.  If $\fu \cap \fk
\not = 0$ then, since $\dim \fu = 1$ we would have $\fu \subset \fk$.
But then $||z||^2$ would be a component of the $K$-moment map $\Phi$
and so $\Phi\inv (0)$ would only contain zero. 

Since $\Phi$ is homogeneous, the level set $\Phi\inv (0)$ is a cone on
$\Phi\inv (0) \cap S^{2n-1}$ where $S^{2n-1}$ is the standard round
sphere in $\C^n$ of radius 1, $S^{2n-1} = \{ z \mid ||z||^2 = 1\} $.
Hence the reduced space $V_0 $ is a cone on the set $L:= (\Phi\inv (0)
\cap S^{2n-1} )/K$, i.e., $V_0 = {\overset \circ c } (L):= (L \times
[0, \infty))/\sim$ where $(x, 0) \sim (x', 0)$ for all $x, x'\in L$.
The vertex $*$ of the cone corresponds to $0\in V$.  Moreover, (see
\cite[Corollary~6.12]{SL}) the set $L$ is a stratified space; it is
the link of the singularity of $V_0$ at $*$.\footnote{Calling $L$ the
link is slightly nonstandard. Strictly speaking we should call $L$ the
link only if the set $\{*\}$ is a stratum; that is, if the set of
fixed points $V^K$ is only the origin.}  The stratifications of $L$
and of $V_0$ are related: given a stratum $\cS$ of $L$, the set $\cS
\times (0,\infty) \subset {\overset \circ c } (L)$ is a stratum of
$V_0$.\footnote{There is one exception: if $V^K\not = \{0\}$ then
$(V^K\cap S^{2n-1})/K = V^K\cap S^{2n-1}$ is a stratum of $L$, but $
(V^K\cap S^{2n-1} )\times (0, \infty) = V^K
\smallsetminus \{0\}$, rather than $V^K$ which is a stratum of $V_0$.
See previous footnote.}

By \cite[Theorem~5.3]{SL} the action of the circle $U$ on $L$ is
locally free and preserves the stratification of $L$.  The quotient
$\cL := L/U$ is again a stratified space.  In fact, $\cL$ is a symplectic
stratified space since it is a reduction of $\C^n$ by the action of
$K\times U$.  The space $\cL$ is called the {\bf symplectic link} of
$*$ in the symplectic stratified space $V_0$.

\begin{remark}
The symplectic link $\cL$ has two natural decompositions.  There is a
stratification of $\cL$ into manifolds as a symplectic stratified
space.  There is also a coarser decomposition: since the action of $U$
on $L$ is locally free and preserves the stratification of $L$, the
quotients of the strata of $L$ form a decomposition of the symplectic
link $\cL$ into symplectic orbifolds: $\cL = \coprod _{\cS \subset L}
\pi(\cS)$ where $\cS$ are strata of $L$ and $\pi: L \to \cL$ is the
$U$-orbit map.  We will use the coarser decomposition.
\end{remark}
\begin{remark}
Note that since the Hamiltonian action of $T$ on $V$ commutes with the
action of $U\times K$, it descends to a Hamiltonian action on the
symplectic link $\cL$.
\end{remark}

Finally recall a description of the symplectic structure on the strata
of the reduced space $V_0$ \cite[Theorem~5.3]{SL}: For each stratum
$\calS $ of $L$ there exists a connection one-form $A_\cS$ on the
Seifert fibration $\calS \to \calS /U $ such that the curvature of
$A_\cS$ is a symplectic form.  The reduced symplectic form on $\calS
\times (0, \infty)$ is $d (s A_\cS)$, where $s$ denotes
the natural coordinate on $(0,\infty)$.  There is no loss of
generaltiy in assuming that the connections $A_\cS$ are $T$-invariant.

We now  study one (connected) stratum $P$ of the link $L$.  Denote
by $B$ the quotient of $P$ by the action of $U$: $B= P/U$.  Then $U
\to P \stackrel{\pi}{\to} B$ is a Seifert fiber bundle.  Denote the
connection one form on $P$ by $A$, so that the symplectic form on
$P\times (0, \infty)$ is $d (sA)$.  By assumption $B$ is a symplectic
orbifold.

\begin{lemma} \label{lemma8}
Let $S^1\to P \stackrel{\pi}{\to} B$ be a Seifert fibration over an
even dimensional closed orbifold $B$.  Suppose there exists a connection
one-form $A$ on $P$ so that the form $d(sA)$ on $P\times (0, \infty)$
is symplectic ($s$ is the coordinate on $(0,\infty)$).  Suppose
further that a torus $T$ acts on $P$ without fixed points, and that the
action of $T$ commutes with the action of $S^1$ and preserves the
connection $A$.  Then the action of $T$ on $(P\times (0, \infty), d
(sA))$ is Hamiltonian.  Let $F$ denotes a corresponding moment map.

Then for a generic vector $Y$ in the Lie algebra $\ft$ of $T$ and for
any $s_0\in (0, \infty)$ the set $\pi \inv (B^T) \times \{s_0\}$
consists of nondegenerate periodic manifolds of the Hamiltonian $H =
\la F, Y\ra$ on the energy surfaces $\{ H = E\}$ for appropriate
$E$'s.
\end{lemma}
\begin{proof}
Since the action of $T$ preserves the one-form $sA$, the action of $T$
on $(P\times (0, \infty), d (sA))$ is Hamiltonian.  The action of
$S^1$ on $(P\times (0,\infty), d (sA))$ is also Hamiltonian: $f(p, s)
= s$ is a corresponding moment map.  Consequently $f\inv (1)/S^1$ is a
symplectic orbifold diffeomorphic to $B$; from now on we identify $B$
and $f\inv (1)/S^1$.

The Hamiltonian action of $T$ on $(P\times (0, \infty), d (sA))$
descends to Hamiltonian action on $B$.  Since $B$ is compact and the
action of $T$ is Hamiltonian, the set of fixed points $B^T$ is
nonempty.  In fact $B^T$ is a disjoint union of connected symplectic
suborbifolds of $B$ (see for example \cite{LTol} for more details on
Hamiltonian group actions on orbifolds).  For any point $x\in
\pi\inv (B^T)$, the $T$ orbit $T\cdot x$ is contained in the $S^1$
orbit $S^1 \cdot x$. Since $T$ acts on $P$ without fixed points we in
fact have that $T\cdot x = S^1 \cdot x$ for any $x\in \pi\inv (B^T)$.
Consequently the union of manifolds $\pi\inv (B^T) \times (0,
\infty) $ consists of periodic manifolds of the Hamiltonian $H$.

It remains to check that for a connected component $\Sigma$ of
$\pi\inv (B^T)$, the manifold $(\Sigma \times (0, \infty)) \cap \{H =
E\}$ is a nondegenerate periodic manifold of $H$.  Now the time $t$ map
of the flow of the Hamiltonian vector field of $H$ on $P\times (0,
\infty)$ is given by $(p, s) \mapsto ((\exp tY)\cdot p, s)$ where
$\exp: \ft \to T$ is the exponential map. 

So let $(p,s)$ be a point in $(\Sigma \times (0, \infty)) \cap \{H =
E\}$.  Then $(p, s)$ is a relative $S^1$ equilibrium of $H$.  Hence
the differential $dH$ at $(p, s)$ is proportional to the differential
of the $S^1$ moment map, which is $ds$.  Hence $T_{(p, s)} \{H =E\} =
T_p P$.  Therefore it's enough to compute the differential at $p$ of the
``Poincar\'e map'' $P\to P$, $q \mapsto \exp (\tau Y) \cdot q$,  where
$\tau$ is the smallest positive number with $\exp (\tau Y)\cdot p = p$.

Since the $T$ orbit of $p$ is a circle, the isotropy group of $p$ is
of the form $\Gamma \times T_2$, where $\Gamma$ is a finite abelian
subgroup of $T$ and $T_2$ is a subtorus of $T$ of codimension one.
Moreover, we can split $T$ as $T= T_1 \times T_2$ where $T_1$ is
isomorphic to $S^1$ and contains $\Gamma$.

Let us next assume, to make the exposition simpler, that $\Gamma$ is
trivial.  Then it follows from the slice theorem that a neighborhood
of $p$ in $P$ is $T$ equivariantly diffeomorphic to a neighborhood of
$(1, 0, 0)$ in $T_1 \times \C^m
\times \C^k$ where $m = \dim \Sigma - 1$. Here $T = T_1 \times T_2$
acts on $T_1 \times \C^m
\times \C^k$ by 
$$
(\lambda, t) \cdot (\mu, w_1, \ldots, w_m, z_1, \ldots z_k) = (\lambda
\mu,  w_1, \ldots, w_m, \chi_1 (t) z_1, \ldots \chi_k (t) z_k), 
$$ 
where $\chi_1, \ldots, \chi_k: T_2 \to U(1)$ are nontrivial
characters of $T_2$.

Let $pr_\alpha: T\to T_\alpha $, $\alpha =1,2$ denote the projections.
Then $pr_1 (\exp(\tau Y)) = 1$.  We claim that for all $r$ between $1$
and $k$, $\chi_r (pr_2 (\exp(\tau Y)))$ is of the form $e^{2\pi i
y_r}$ where $y_r$ are {\bf irrational} numbers.  Note that the claim
implies immediately that the algebraic multiplicity of the eigenvalue
1 of the differential of the ``Poincar\'e map'' $q \mapsto \exp (\tau
Y) \cdot q$ is $\dim \Sigma$, hence that $(\Sigma \times (0, \infty))
\cap \{H = E\}$ is a nondegenerate periodic manifold of $H$.

The claim holds because the one-parameter subgroup $\{ \exp (tY) \mid
t\in \R\}$ is dense in $T$.  More specifically, let $e_1, \ldots e_s$
be a basis of the integral lattice of the torus $T$ which is
compatible with the splitting $T= T_1 \times T_2$  (so that  $\{\exp
(te_1) \mid t \in \R \} = T_1$  and $e_2, \ldots, e_s$ is a basis of
the integral lattice of $T_2$).  Then $Y = a_1e_1 + \sum_{j=2}^s a_j
e_j$ for some $a_j\in \R$.  Moreover, since the one-parameter subgroup
defined by $Y$ is dense in $T$, the sum $\sum_{j=1}^s q_j a_j$ is not
a rational number for any $s$ tuple of rational numbers $(q_1, \ldots,
q_s)$.    Since $\exp (pr_1 (\tau Y)) = 1$, $a_1 = \pm
\frac{1}{\tau}$.  Consequently 
$$
\chi_r (pr_2 (\exp(\tau Y))) = \chi_r
(\pm \sum_{j=2}^{s} \frac{a_j}{a_1} e_j) = e^{2\pi i (\pm
\sum_{j=2}^{s}\frac{d\chi_r (e_j)}{ a_1} a_j)}
$$ 
and the claim follows (note that $d\chi_r (e_j)$ are integers).

If the group $\Gamma$ is not trivial, then a neighborhood of $p$ in
$P$ is modeled by the quotient $(T_1 \times \C^m \times
\C^k)/\Gamma$, where $\Gamma$ acts on $T^1$ by multiplication and on
$\C^m \times \C^k$ linearly by $m+k$ characters, so that the actions
of $T$ and $\Gamma$ commute.  The same argument as above still works:
for the Poincar\'e map on the quotient to have an eigenvector in
$\C^k$ with eigenvalue 1 it is necessary for $\chi_r (pr_2 (\exp(\tau
Y)))$ to be a root of unity.  But this is impossible as we have seen.
This proves Lemma~\ref{lemma8}.
\end{proof}

In fact in proving Lemma~\ref{lemma8}, we have proved more:
\begin{proposition}\label{prop5}
Let $(V, \omega)$ be a symplectic vector space with a linear action of
a compact Lie group $K$ and a corresponding homogeneous moment map
$\Phi: V \to \fk^*$.  Let $q\in C^\infty (V)^K$ be a $K$ invariant
positive definite quadratic Hamiltonian; let $q_0$ be the
corresponding reduced Hamiltonian on $V_0 = \Phi\inv (0)/K$.  Let
$T\subset \SP (V, \omega)$ be the torus generated by $q$.  Let $L$ be
the link of the point $*\in V_0$ corresponding to $0$, let $\cL$ be
the symplectic link and let $\pi: L \to \cL$ be the orbit map.

Then for every stratum $\cS$ of the link $L$ such that the fixed point
set $\pi(\cS)^T$ is nonempty and for every $E>0$ the manifold 
$$
\{ q_0 = E\} \cap (\cS\times (0, \infty))
$$
contains a weakly nondegenerate periodic manifold $C$ of $q_0$: $C =
\pi\inv (\pi(\cS)^T)$.  If the orbifold $\pi(\cS)^T$ is compact, then the
manifold $C$ is compact as well.
\end{proposition}

Note that by construction the circle action on the periodic manifolds
$C$ in Proposition~\ref{prop5} is simply the action of the circle $U$.
Hence $C/S^1 = \pi\inv (\pi(\cS)^T)/U =\pi(\cS)^T$.  It follows from a
result of Weinstein \cite[p. 247]{W2} that in
Theorem~\ref{theorem_linear} the number of periodic orbits of the
Hamiltonian $h_0$ on a given energy surface $\{h_0 = E\}$ is bounded
below by
\begin{equation}\label{eqN1}
N_1= \sum \text{Cat} (\pi (\cS)^T)
\end{equation}
where the sum is taken over all strata $\cS$ of the link $L$ such that
the sets $\pi (\cS)^T$ are {\bf compact}.  Since the link $L$ is
compact, the closed strata of $L$ must be compact.  It follows that
the number $N_1$ in equation (\ref{eqN1}) is positive.

The bound given by (\ref{eqN1}) is somewhat unsatisfactory --- it
ultimately depends on the Hamiltonian, while no such dependence is
present in Weinstein's nonlinear normal modes theorem
(Theorem~\ref{theorem_W} above).  We will see in Lemma~\ref{lemma_LS}
below that $\text{Cat} (\pi (\cS)^T) \geq \text{Cat} (\pi (\cS))$ for
any closed stratum $\cS$ of the link $L$.  Consequently $$ N_1 \geq
\sum \text{Cat} (\pi (\cS)), $$ where the sum is taken over all closed
strata $\cS$ of the link $L$.  This will finish our proof of
Theorem~\ref{theorem_linear} hence of Theorem~\ref{main_theorem}.

\begin{lemma}\label{lemma_LS}
Let $B$ be a closed symplectic orbifold with a Hamiltonian action of a
torus $T$.  Then the Liusternik-Schnirelman category $\text{Cat}
(B^T)$ of the set of $T$-fixed points is bounded below by the
Liusternik-Schnirelman category of $B$: $\text{Cat} (B^T) \geq
\text{Cat}(B)$.
\end{lemma}
\begin{proof} We use open sets in our  definition of
the category.  Let $f: B \to \R$ be a generic component of the moment
map for the action of $T$ on $B$ so that $B^T$ is precisely the set of
critical points of $f$.  The function $f$ is Bott-Morse.  Therefore
$B$ decomposes as a disjoint union of the unstable orbifolds $W_1,
\ldots, W_k$ of the gradient flow of $f$.  Clearly $\text{Cat}(\coprod
W_k) = \text{Cat} (B^T)$.  Now ``thicken'' $W_j$'s by replacing them
with their tubular neighborhoods $\tilde{W}_j$ inside $B$.  Then
$\text{Cat} (\tilde{W}_j) = \text{Cat} (W_j)$, the sets
$\tilde{W}_j$'s are open and $\cup_j \tilde{W}_j = B$.  Hence
$\text{Cat} (B) \leq \sum \text{Cat} (\tilde{W}_j) = \text{Cat}(B^T)$.
\end{proof}

We end the paper by describing a practicable method for checking the
existence of a symplectic slice $\Sigma$ through a relative
equilibrium $m$ of a symmetric Hamiltonian system $(M,\omega _M , \Phi:
M \to \fg^*, h \in C^\infty (M)^G)$ so that the Hessian $d^2
(h|_\Sigma) (m)$ is positive definite.

Recall that if a point $m$ is a relative equilibrium of a symmetric
Hamiltonian system, then there exists a vector $\eta \in \fg$ so that
\begin{equation}\label{eq_rel_eq}
d (h - \la \Phi, \eta \ra) (m) = 0. 
\end{equation}
Then the Hessian $d^2 (h - \la \Phi, \eta \ra) (m)$ is a well-defined
quadratic form, which we will use shortly.
The vector $\eta$ is not unique: for every $\zeta$ in the Lie algebra
$\fg_m$ of the isotropy group of $m$, the vector $\eta + \zeta$ also
satisfies  $d (h - \la \Phi, \eta + \zeta \ra) (m) =
0$.  It is not hard to show that $\eta $ has to lie in the isotropy
Lie algebra $\fg_\mu$ where $\mu = \Phi (m)$.

Since by assumption the action of $G$ on $M$ is proper, the
isotropy group $G_m$ is compact. Hence we can choose a $G_m$ invariant
inner product on the Lie algebra $\fg$ and use it to define an
orthogonal complement $\fm$ of $\fg_m$ in $\fg_\mu$.  There exists a {\em
unique} vector $\eta \in \fm$ so that (\ref{eq_rel_eq}) holds.  The
vector $\eta$ is called the {\bf orthogonal velocity} of the relative
equilibrium $m$ \cite{OR1}.

\begin{proposition} \label{last prop}
Let $m$ be a relative equilibrium of a symmetric Hamiltonian system
$(M,\omega _M , \Phi: M \to \fg^*, h \in C^\infty (M)^G)$ and let
$\eta\in \fg_\mu$ be the orthogonal velocity of $m$ with respect to
some $G_m$ invariant inner product on $\fg$ (where $\mu = \Phi (m)$).
Suppose that the quadratic form $\left. d^2 (h - \la \Phi, \eta
\ra)(m) \right|_{\ker d\Phi (m)}$ is semi-definite of maximal possible rank
(the dimension of a symplectic slice at $m$).  Then there exists a
symplectic slice $\Sigma$ through $m$ such that the form $d^2
(h|_\Sigma)(m)$ is positive definite.
\end{proposition}
\begin{proof}
The proof is a standard computation that uses the local normal form of
the moment map of Marle and of Guillemin and Sternberg
\cite{marle:model, GS}.  Similar computations are carried out in
\cite[p. 1643]{LS} and in \cite{OR1}.
We use the version of the normal form theorem described in
\cite[pp. 211--215]{BL} which we now recall without proofs:\\

Let the symbols $(M, \omega)$, $G$, $\Phi:M \to \fg^*$, 
$\mu = \Phi (m)$,  $\fg_m$, $\fg_\mu$ and $\fm$ have the same meaning as above.

The null directions of the restriction $\omega (m)|_{\ker d\Phi (m)}$
is $T_m (G_\mu \cdot m)$.  Hence  $V= \ker d\Phi (m)/ T_m (G_\mu \cdot
m)$ is naturally a symplectic vector space.  Denote the corresponding
symplectic form by $\omega _V$.  Moreover, the linear action of the isotropy
group $G_m$ on $\ker d\Phi (m)$ descends to a linear symplectic action
on $(V, \omega _V)$.  Denote the corresponding homogeneous moment map
by $\Phi _V$.

The $G_m$ invariant inner product on $\fg$ chosen above defines $G_m$
equivariant embeddings: $i: \fg_m^* \to \fg^*$ and $j:\fm ^* \to
\fg^*$.  Note that by construction of $i$ and $\fm$ we have that $\la
i(\ell), \eta\ra = 0$ for any $\ell \in \fg_m^*$ and any $\eta \in \fm$.

There exists a closed two-from $\sigma$ on the associated bundle $Y=
G\times _{G_m} (\fm^* \times V)$ and an open $G$ equivariant embedding
$\psi$ of a neighborhood the zero section of $Y\to G/G_m$ into $M$
with the following properties.
\begin{enumerate}
\item  $\psi ([1, 0, 0]) = m$.

\item $\psi ^* \omega  = \sigma$.

\item  $(\psi ^* \Phi)([g, \lambda, v]) = Ad^\dagger (g)(\mu + j (\lambda)
+ i(\Phi _V (v)))$ for all $[g, \lambda, v]\in G\times _{G_m} (\fm^*
\times V)$. 

\item  The embedding $\iota: (V, \omega _V) \hookrightarrow (G\times
_{G_m} (\fm^* \times V), \sigma)$, $\iota(v) = [1, 0, v]$ is
symplectic.  Consequently for a small enough neighborhood $U$ of 0 in
$V$, $ \psi (\iota (U))$ is a symplectic slice through $m$.\\

\end{enumerate}

We now prove that $\Sigma = \psi (\iota (U))$ is the desired
symplectic slice.  Since $d(h - \la \Phi, \eta\ra)(m) =0$, the Hessian
$d^2(h - \la \Phi, \eta\ra)(m) $ is well-defined and behaves well
under restrictions.  In particular, $d^2(h - \la \Phi,
\eta\ra)(m)|_{T_m (G_\mu \cdot m)} = d^2 (h - \la \Phi,
\eta\ra)|_{G_\mu \cdot m} ) (m)$.  Since $h$ is $G$ invariant and
since for any $a\in G_\mu$ we have $\la \Phi, \eta \ra (a\cdot m) =
\la Ad^\dagger (a) \Phi (m), \eta\ra = \la Ad^\dagger (a) \mu, \eta\ra
= \la \mu, \eta \ra = \la \Phi, \eta \ra (m)$. Therefore $(h - \la \Phi,
\eta\ra)|_{G_\mu \cdot m} $ is constant and hence $d^2(h - \la \Phi,
\eta\ra)(m)|_{T_m (G_\mu \cdot m)} = 0$.

Since the null directions of $\omega (m)|_{\ker d\Phi (m)}$
is $T_m (G_\mu \cdot m)$, it follows that for any symplectic slice
$\Sigma'$ through $m$ which is tangent to $\ker d\Phi (m)$, we have
$$
T_m \Sigma' \oplus T_m (G_\mu \cdot m) = \ker d\Phi (m).
$$
Combining this with the previous computation we see that the rank of
$\left. d^2 (h - \la \Phi, \eta \ra)(m) \right|_{\ker d\Phi (m)}$ is
at most $\dim \Sigma'$.  Thus by assumption $ 
\left. d^2 (h - \la \Phi, \eta \ra)(m) \right|_{T_m \Sigma'}$ is positive
definite for any symplectic slice which is tangent to $\ker d\Phi
(m)$.  It is easy to check that the manifold $ \psi (\iota (U))$ is
such a slice.   We finally show that 
$$
\left.d^2 (h - \la \Phi, \eta \ra)(m) \right|_{ \psi (\iota (U))} = d^2
(h|_{\psi (\iota (U))}) (m).
$$
Now for any $u\in U$ we have $\la \Phi,
\eta (\psi (\iota (u))\ra = \la \psi^* \Phi, \eta\ra ([1, 0, u])= \la \mu
+ j(0) + i(\Phi_V (u)), \eta\ra = \la \mu, \eta\ra$ (since $\la
i(\fg_m^*), \eta\ra = 0$ by construction of $i$ and $\eta$).  
Therefore $
\left. d^2 (h - \la \Phi, \eta \ra)(m) \right|_{ \psi (\iota (U))}= d^2 ((h -
\la \Phi, \eta\ra)|_{ \psi (\iota (U))}) (m) = d^2 ( h|_{ \psi (\iota
(U))} - \la \mu, \eta\ra) (m) =d^2 (h|_{\psi (\iota (U))}) (m)$.
\end{proof}

\end{document}